\documentclass[12pt]{article}

\usepackage[margin=1in]{geometry} 
\usepackage{amsmath,amsthm,amssymb,graphicx,mathtools,tikz,hyperref}
\usepackage[bottom]{footmisc}
\usetikzlibrary{positioning}

\newcommand{\footremember}[2]{%
    \footnote{#2}
    \newcounter{#1}
    \setcounter{#1}{\value{footnote}}%
}
\newcommand{\footrecall}[1]{%
    \footnotemark[\value{#1}]%
}

 \hypersetup{
 colorlinks,
 linkcolor=blue
 }

\begin{document}
\date{\today}

\title{Derivation of Analytic Formulas for the Sample Moments of the Sample Correlation over Permutations of Data}
%\author{Marc Jaffrey$^{1,3}$ \\
%Michael Dushkoff$^{2,3}$}
\author{
  Marc Jaffrey\footremember{uw}{University of Washington (\href{mailto:mjaffrey@uw.edu}{mjaffrey@uw.edu})}\footremember{rf}{RootFault LLC} \and
  Michael Dushkoff\footremember{rit}{Rochester Institute of Technology (\href{mailto:mad1841@rit.edu}{mad1841@rit.edu})}\footrecall{rf}
}

\maketitle

% % % % % % % % % % % % % % % % % % % % % % % % %
% Abstract
\begin{abstract}
\noindent Pearson's correlation $\rho$ is among the mostly widely reported measures of association. The strength of the statistical evidence for linear association is determined by the p-value of a hypothesis test. If the true distribution of a dataset is bivariate normal, then under specific data transformations a t-statistic returns the  exact p-value, otherwise it is an approximation. Alternatively, the p-value can be estimated by analyzing the distribution of the sample correlation, $r$, under permutations of the data. Moment approximations of this distribution are not as widely used since estimation of the moments themselves are numerically intensive with greater uncertainties.  In this paper we derive an inductive formula allowing for the analytic expression of the sample moments of the sample correlation under permutations of the data in terms of the central moments of the data. These formulas placed in a proper statistical framework could open up the possibility of new estimation methods for computing the p-value.
\end{abstract}

%%%%%%%%%%%%%%%%%%%%%%%%%%%%%%%%%%%%%%%%%%%%%%
% Introduction
%%%%%%%%%%%%%%%%%%%%%%%%%%%%%%%%%%%%%%%%%%%%%%
\section{Introduction}
\begin{figure}
    \begin{center}
      \advance\leftskip-3cm
      \advance\rightskip-2.3cm
      \includegraphics[scale=0.25]{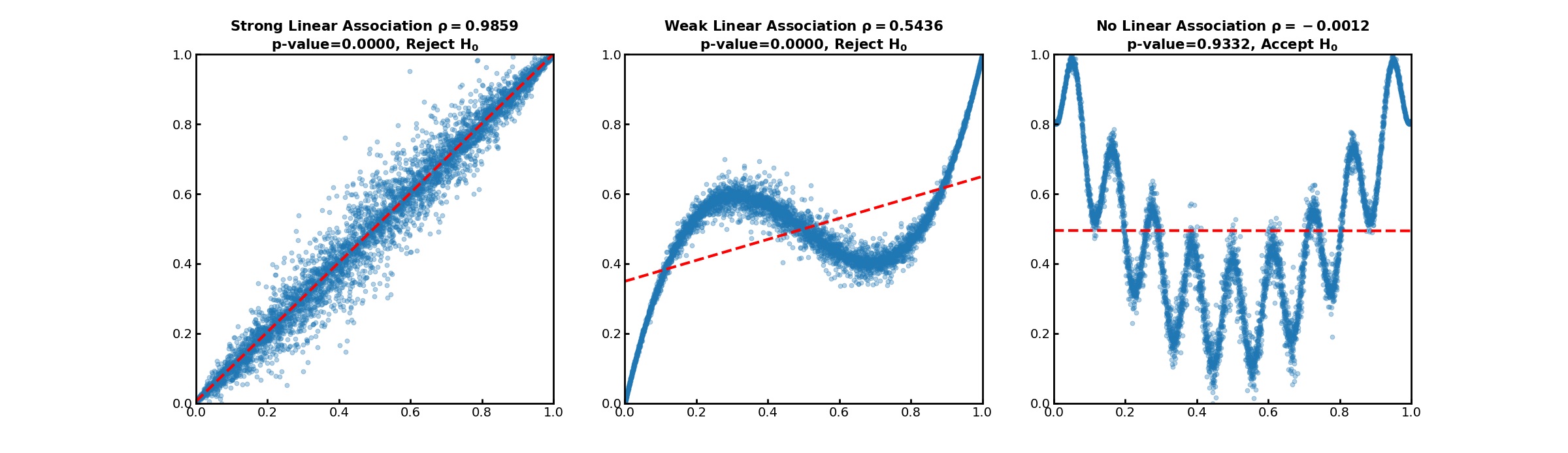}
      \caption{ All three synthetic examples clearly exhibit strong functional relationships. Under a hypothesis test for the presence of a linear relationship under $\rho$, the first two examples strongly reject the null hypothesis, both with p-values $< \mathit{1 \times 10^{-4}}$, while in the last example the null hypothesis is accepted. Even though the middle example reports a  moderate linear association, it is clearly non-linear being a cubic polynomial. Pearson's correlation is perhaps better interpreted as a measure of the suitability of a linear approximation to the data.}
      \label{fig:pc}
  \end{center}
\end{figure}

Pearson's product-moment correlation \cite{pearson1895correlation},
\begin{equation}
    \rho = \frac{COV(X,Y)}{\sigma_X\sigma_Y}
    \label{ppmc}
\end{equation}
the generalization of Galton's regression coefficient,  \cite{galton1886regression,stigler1989francis} is perhaps the most important measure of data association, orienting scientists in the direction of discovery since its introduction.  As an inferred measure of \textit{linear} association between two variables, it is a primary tool in statistical analysis and the search for meaningful variable relationships, (see Figure \ref{fig:pc}). The subject of its interpretation and nuances surrounding its use have been extensively covered in the literature with still plenty of room for debate today, \cite{lee1988thirteen,puth2014effective}.

Given a dataset with two variables, in hypothesis testing,
\begin{align}
    & H_0: \rho=0 \\
    & H_1: \rho \neq 0 \nonumber
\end{align}
the p-value of the test determines the strength of the statistical evidence for whether to accept or reject $H_0$ \cite{pitman1937significance}. Widely contentious, the use of p-value for determining statistical significance is under critical debate, a subject that is extensively addressed in the recent literature \cite{chawla2017p,head2015extent,held2018p,krueger2017heuristic,krueger2019putting,kuffner2019p}. Notwithstanding this controversy, p-value is still in widespread use today, and as such, methods of its computation are relevant to discuss. 
\\
\\
Presented is a branching inductive formula for computing the sample moments of the sample correlation coefficient $r$, 
\begin{equation}
    r = \frac{\sum_i (x_i -\mu_x)(y_i-\mu_y)}{ \sqrt{\sum_k(x_k-\mu_x)^2\sum_j(y_j-\mu_y)^2}},
    \label{ssc}
\end{equation}
over permutations of the data in terms of the central moments of the data, denoted
\begin{equation}
    \langle r_{\Pi}^k\rangle = \frac{1}{|\Pi|} \sum_{\pi \in \Pi} r_{\pi}^{k}
    \label{mdef}
\end{equation}
for $k \in \mathbb{N}$, where $\Pi$ is the set of all possible permutations on the data and $r_{\pi}$ denotes the sample correlation over the permutation $\pi \in \Pi$ of the data. These formulas are interesting by themselves and placed in a proper statistical framework could potentially open the door to the possibility of computationally efficient methods for computing the p-value of a hypothesis test of $\rho$.
%These exact formulas, though interesting by themselves, open the door to the possibility of computationally efficient methods for evaluating the p-value of a hypothesis test of $\rho$ under the permutation distribution for the sample correlation using moment approximations.

% https://www.mathworks.com/matlabcentral/fileexchange/727-exportfig

%%%%%%%%%%%%%%%%%%%%%%%%%%%%%%%%%%%%%%%%%%%%%%
% 
%%%%%%%%%%%%%%%%%%%%%%%%%%%%%%%%%%%%%%%%%%%%%%
\section{Main Result: Inductive Formula for $\langle r_{\Pi}^k\rangle$}
\subsection{Notation} Given a dataset $D_n = \lbrace ( x_i,y_i )|\ i=1,..,n \rbrace \subset \mathbb{R}^2$, let $\Pi=Perm(n)$ be the set of all permutations on $\{ 1,\dots,n\}$. For $\pi=(\pi_1,\dots,\pi_n)\in \Pi$, defined such that $i\underset{\pi}{\to} \pi_i$, define
\begin{equation}
 \pi(D_n)=\lbrace(x_i,y_{\pi_i})|i=1,..,n\rbrace, 
 \label{py}
\end{equation}
where the permutation $\pi$ acts on the $y$-coordinate alone. Denote Pearson's sample correlation over a permutation $\pi$ of the data as $r_\pi = r(\pi(D_n))$
and define 
\begin{eqnarray}
   r_{\Pi}: [-1,1] \to [0,1],
\end{eqnarray}
the distribution of $r$ over the permutations of $D_n$.
Lastly, let $\hat{x}_i \equiv (x_i -\mu_x)$ and $\hat{y}_i \equiv   (y_i - \mu_y)$.

\subsection{Main Result}
Given $D_n$, fix $k$. The $k^{th}$ moment of $r_{\Pi},$  is given by
\begin{equation}
    \langle r_{\Pi}^k \rangle   =   \frac{1}{n^k n!\hat{\sigma}_x^k \hat{\sigma}_y^k} \sum_{m=1}^{k} \ \ \sum_{n_1+...+n_m=k} \binom{k}{n_1,..,n_m}^* X_{(n_1,..,n_m)}^{m,k} (n-m)!Y_{(n_1,..,n_m)}^{m,k}\cdot h_{n,m}
    \label{pcpm}
\end{equation}
$X_{(n_1,..,n_m)}^{m,k} \text{ and } Y_{(n_1,..,n_m)}^{m,k} $ are branching inductive formulas, such that for $1<m\leq k$ and for fixed non-zero positive integers $n_1+...+n_m=k$ we have:
\begin{equation}
    X_{(n_1,..,n_m)}^{m,k}  =   \ n \langle  \hat{x}^{n_m} \rangle X_{(n_1,...,n_{m-1})}^{m-1,k-n_m}  - \sum_j X_{(n_1 + n_m\delta_{ij},...,n_{m-1}+n_m\delta_{(m-1)j})}^{m-1,k} 
    \label{inducx}
\end{equation}
where $\delta_{ij}$ is the standard delta function, and similarly defined
\begin{equation}
    Y_{(n_1,..,n_m)}^{m,k}  =   \ n \langle  \hat{y}^{n_m} \rangle Y_{(n_1,...,n_{m-1})}^{m-1,k-n_m}  - \sum_j Y_{(n_1 + n_m\delta_{ij},...,n_{m-1}+n_m\delta_{(m-1)j})}^{m-1,k}
    \label{inducy}
\end{equation}
For $m=1$, we have
\begin{equation}
X_{(k)}^{1,k}  =  n\langle \hat{x}^k \rangle  \text{ and }
     Y_{(k)}^{1,k} =  n\langle \hat{y}^k \rangle
\end{equation}
The stared multinomial coefficient
%\hfill\\
\begin{equation}
    \binom{k}{n_1,\dots,n_m}^* = 
\frac{1}{d_1!\cdots d_r!}\binom{k}{n_1,\dots,n_m}
\label{multi}
\end{equation}
is an adjustment of the usual multinomial coefficient accounting for degeneracy in $n_1,...,n_m$. Partitioning $n_1,...m_m$ into subsets, $g_1,..,g_r$, by the equivalence relation $n_i\equiv n_j \iff n_i=n_j$, then $d_i=|g_i|$. This degeneracy leads to an over counting represented by the multinomial coefficient which is correct by dividing out $d_1!\cdots d_r!$.

The term $h_{n,m}$ accounts for the inability to compute higher order terms in the sum when the number of data points is less than the moment order being computed,  by setting them to zero in the formula,
\begin{equation}
h_{n,m} = \left\{ \begin{array}{rcl}
 0 & \mbox{for}
& n-m < 0 \\ 1 & \mbox{for} & n-m \geq 0
\end{array}\right. 
\end{equation}
Lastly, for notational simplicity, in \eqref{pcpm} the following convention is employed,
\begin{equation}
    \hat{\sigma}_z = \sqrt{\frac{1}{n}\sum (z_i -\mu_z)^2}
    \label{sigma}
\end{equation}

% % % % % % % % % % % % % % % % % % % % % % % % 
% Derivation outline
\section {Derivation}
Starting from \eqref{mdef} and \eqref{ppmc}, we have
\begin{eqnarray}
    \langle r_{\Pi}^k \rangle  & = & \frac{1}{n!}\sum_{\pi \in \Pi} r_{\pi}^k\\
    & = &\frac{1}{n!}\sum_{\pi\in \Pi}\Bigg[  \frac{\sum_i (x_i -\mu_x)(y_{\pi_i}-\mu_y)}{
    \sqrt{\sum_k(x_k-\mu_x)^2\sum_j(y_{\pi_j}-\mu_y)^2}}\Bigg ]^k 
    \label{pp0}
    \\
    &  = & \frac{1}{n^k n!\hat{\sigma}_x^k \hat{\sigma}_y^k} \sum_{\pi\in \Pi} \sum_{i_i,..,i_k} \hat{x}_{i_1}...\hat{x}_{i_k} \hat{y}_{\pi_{i_1}}...\hat{y}_{\pi_{i_k}} 
    \label{pp1}
   % \\
%    & = &\frac{1}{n^k n!\sigma_x^k \sigma_y^k}  \sum_{i_i,..,i_k} \hat{x}_{i_1}...\hat{x}_{i_k} \sum_{\pi\in \Pi} \hat{y}_{\pi_{i_1}}...\hat{y}_{\pi_{i_k}} 
 %   \label{pp2}
\end{eqnarray}
Reorganizing the summation over the indices, let $1\leq m \leq k$ be the number of distinct indices and the non-zero positive integers $n_1+\cdots+n_m=k$ their multiplicity, then \eqref{pp1} becomes
\begin{eqnarray}
      \frac{1}{n^k n!\hat{\sigma}_x^k \hat{\sigma}_y^k}  \sum_{m=1}^{k}\   \sum_{n_1+...+n_m=k} \ \sum_{i_1\neq ... \neq i_m} \binom{k}{n_1,..,n_m}^* \hat{x}_{i_1}^{n_1}...\hat{x}_{i_m}^{n_m} \sum_{\pi} \hat{y}_{\pi_{i_1}}^{n_1}...\hat{y}_{\pi_{i_m}}^{n_m}.
     \label{pp3} 
\end{eqnarray}
The adjustment to the multinomial coefficient is addressed later. 

%%%%%%%%%%%%%%%%%%%%%%%%%%%%%%%%%%%%%%%%%%%%%%%%%%%%%
\subsection{Summation over the Set of Permutations}
Beginning with the summation over the permutations, fix $m$, $i_1\neq\cdots\neq 1_m$,  and  $n_1+...+n_m=k$,
\begin{equation}
    \sum_{\pi \in \Pi} \hat{y}_{\pi_{i_1}}^{n_1}...\hat{y}_{\pi_{i_m}}^{n_m}
    \label{permysum}
\end{equation}
This sum is independent of the choice of indices $i_1\neq...\neq i_m$ inherited from the summation over the $x$-coordinates in \eqref{pp3}. To see this, take any two subsets of indices $i_i\neq...\neq i_m$ and $j_1\neq...\neq j_m$ and let $\sigma_{ij}\in \Pi$ be any permutation satisfying $\sigma_{ij}(i_l) = j_l$ for $l= 1,..,m$,  ($\sigma_{ij}$ is not unique). Taking $\Pi$ as a group \cite{armstrong2013groups}, it follows that $\Pi\circ\sigma_{ij} = \Pi$ and consequently,
\begin{eqnarray}
    \sum_{\pi \in \Pi} \hat{y}_{\pi_{i_1}}^{n_1}...\hat{y}_{\pi_{i_m}}^{n_m} & = & \nonumber \\
    & = & \sum_{\pi\circ\sigma_{ij}\in \Pi} \hat{y}_{(\pi \circ\sigma_{ij})_{i_1}}^{n_1} \cdots\hat{y}_{(\pi\circ\sigma_{ij})_{i_m}}^{n_m}\\
    & = & \sum_{\pi \in \Pi} \hat{y}_{\pi_{j_1}}^{n_1}...\hat{y}_{\pi_{j_m}}^{n_m} \nonumber
\end{eqnarray}
Thus we can shift the summation to the first $m$ indices, which we denoted as,
\begin{equation}
     (n-m)!Y_{(n_1,...,n_m)}^{m,k} =\sum_{\pi\in\Pi} \hat{y}_{\pi_1}^{n_1}...\hat{y}_{\pi_{m}}^{n_m}
     \label{sump}
\end{equation} 
The necessity of $(n-m)!$ will be come clear in a moment.
Next, place the $N=n!$ elements of $\Pi$ into the rows of a matrix and define the operation of $\Pi$ on the $y$ coordinates as
 \begin{eqnarray}
 \Pi(\hat{y})
 = \begin{pmatrix}
 \hat{y}_{\pi_1^1} & \cdots \cdots & \hat{y}_{\pi_n^1} \\
 \vdots  & \ddots & \vdots\\
  \hat{y}_{\pi_1^N} & \cdots \cdots & \hat{y}_{ \pi_n^N}
 \end{pmatrix},
 \label{matrixp}
 \end{eqnarray}
 \eqref{sump} is now interpreted as the sum over the rows of the product of the element in the first $m$ columns raised to the appropriate power. As the order of the rows does not matter for the summation, rearrange \eqref{matrixp} into the canonical form,
\begin{equation}
 \Pi^*(\hat{y}) =
  \begin{pmatrix}
    \hat{y_1} &
    \boxed{\begin{minipage}[c][1.85cm][b]{9cm}\centering
      $$\begin{matrix}
         \hat{y_2} & \boxed{\begin{minipage}[c][1.1cm][b]{8cm}\centering
          $$\begin{matrix}
             \hat{y_3} & \boxed{\begin{minipage}{7.1cm}\centering $\cdots$ \end{minipage}}\\
            \vdots\\
             \hat{y_n} &
          \end{matrix}$$
        \end{minipage}}\\
        \vdots
        \\
         \hat{y_n} & \boxed{\begin{minipage}[c][0.35cm][b]{8cm}\centering $\vdots$ \ \
          \boxed{\begin{minipage}{7.1cm}\centering$\cdots$\end{minipage}}
        \end{minipage}}
      \end{matrix}$$
    \end{minipage}}
    \\
    \\
     \hat{y_2} &
    \boxed{\begin{minipage}[c][1.6cm][b]{9cm}\centering
      $$\begin{matrix}
       \hat{y_1} &
        \boxed{\begin{minipage}[c][0.35cm][b]{8cm}\centering
          $\vdots$ \ \ 
          \boxed{\begin{minipage}{7.1cm}\centering$\cdots$\end{minipage}}
        \end{minipage}}\\
        \\
         \hat{y_3} & \boxed{\begin{minipage}{8cm}\centering
          $\cdots$
        \end{minipage}}\\
        $\vdots$
        \\
        \hat{y_n} & \boxed{\begin{minipage}{8cm}\centering
          $\cdots$
        \end{minipage}}
      \end{matrix}$$
    \end{minipage}}
    \\
    \vdots
    \\
    \\
     \hat{y_n} &
    \boxed{\begin{minipage}{9cm}\centering
      $\cdots$
    \end{minipage}}
 \end{pmatrix}
 \label{CanM}
\end{equation}
The key observation is that for every distinct sequence of points $\{y_{p_1},...,y_{p_{j}}\}$ in the first $j-1$ columns, associated with it in the $j^{th}$ column are the remaining $\{ y_1,...,y_n\} \backslash  \{y_{p_1},...,y_{p_{j}}\}$ each with multiplicity $(j-1)!$, arising from the remaining $j-1$ columns ignored. Out of this observation the branching inductive formula is built.

For the case $m=1$, summing along the rows of the canonical form, and a slight abuse of notation, 
\begin{equation}
     (n-1)! Y_{(k)}^{1,k} = \sum_{\pi \in\Pi} \hat{y}_{\pi_{1}}^{k} =
     \sum_{Col_1} 
\begin{pmatrix}
     \hat{y_1}^k &
    \boxed{\begin{minipage}{4cm}\centering
      $\vphantom{\int^0}\smash[t]{\vdots}$
    \end{minipage}}\\
    \\
         \hat{y_2}^k &
    \boxed{\begin{minipage}{4cm}\centering
      $\vphantom{\int^0}\smash[t]{\vdots}$
    \end{minipage}}\\
    \vdots \\
     \hat{y_n}^k &
    \boxed{\begin{minipage}{4cm}\centering
      $\vphantom{\int^0}\smash[t]{\vdots}$
    \end{minipage}}
 \end{pmatrix}
\end{equation}
where each factor $\hat{y}_i^{k}$ is repeated $(n-1)!$ times, such that 
\begin{eqnarray}
    \sum_{\pi\in\Pi} \hat{y}_{\pi_{1}}^{k}  & = & \sum_{i=1}^{n} \hat{y}_i^{k}(n-1)!\nonumber \\
    & = & n\cdot \langle \hat{y}^{k} \rangle (n-1)!
    \label{sumi}
\end{eqnarray}
We now proceed to the branching inductive step. Let $1<m\leq k$ and take non-zero positive integers $n_1+..+n_m =k$. We need the following observation, 
\begin{equation}
  \sum_{i \notin \{i_1,...,i_{m-1}\}}\hat{y}_i ^{n_m}\\
  %\sum_{i \neq i_1,...,i_{m-1}}\hat{y}_i ^{n_m}\\
     = \Big [ n\cdot \langle \hat{y}^{n_m}\rangle  - \hat{y}_{i_1}^{n_m} - \hdots -\hat{y}_{i_{m-1}}^{n_m}\Big].
     \label{indstep}
\end{equation}
Returning to \eqref{sump} and reorganizing the sum over distinct sequences of points in the first $m-1$ columns,
\begin{equation}
    (n-m)!Y_{(n_1,...,n_m)}^{m,k}  =  \sum_{\pi} \hat{y}_{\pi_{1}}^{n_1}\cdots\hat{y}_{\pi_{m}}^{n_{m}},
    \label{yh1}
\end{equation}
we obtain
\begin{equation}
      (n-m)!Y_{(n_1,...,n_m)}^{m,k}    =  \sum_{i_1,..,i_{m-1}} \hat{y}_{i_1}^{n_1}\cdots\hat{y}_{i_{m-1}}^{n_{m-1}}\sum_{y_m\neq y_{i_1},...,y_{i_{m-1}}}  \hat{y}_m^{n_m} (n-m)!
    \label{branchmid}
\end{equation}
which by \eqref{indstep} equals
\begin{equation}
     =  \sum_{i_1,..,i_{m-1}} \hat{y}_{i_1}^{n_1}\cdots\hat{y}_{i_{m-1}}^{n_{m-1}}[n\langle \hat{y}^{n_m} \rangle -  \hat{y}_{i_1}^{n_m} \cdots - \hat{y}_{i_{m-1}}^{n_m} ] (n-m)! %\nonumber
    \label{branch}
\end{equation}
From \eqref{branch} the branching inductive formula is clear,
\begin{equation}
    Y_{(n_1,..,n_m)}^{m,k}  =   \ n \langle  \hat{y}^{n_m} \rangle Y_{(n_1,...,n_{m-1})}^{m-1,k-n_m}  - \sum_j Y_{(n_1 + n_m\delta_{ij},...,n_{m-1}+n_m\delta_{(m-1)j})}^{m-1,k}
    \label{inducy2}
\end{equation}
The rationale for separating out the factor $(n-m)!$ in \eqref{sump} is now clear, as the induction step applies only to the sum over indices and not multiplicity arsing from the set of permutations. The multiplicity that arises in the initial sum over permutations is not present in the induction step once the initial sum is collapsed to a sum over indices alone.

%%%%%%%%%%%%%%%%%%%%%%%%%%%%%%%%%%%%%%%%%%%%%%%%%%%%
\subsection{Summation over the $\hat{x}$ terms}
Moving to the summation over the $x$-coordinates,
\begin{equation}
\sum_{i_i,..,i_k} \hat{x}_{i_1}...\hat{x}_{i_k} =
       \sum_{m=1}^{k}  \sum_{n_1+...+n_m=k} \ \  \sum_{i_1\neq ... \neq i_m} \binom{k}{n_1,..,n_m}^* \hat{x}_{i_1}^{n_1}...\hat{x}_{i_m}^{n_m},
       \label{xterms}
\end{equation}
and focusing on 
\begin{equation}
    \sum_{i_1\neq ... \neq i_m} \hat{x}_{i_1}^{n_1}...\hat{x}_{i_m}^{n_m},
\end{equation}
this is the same exact sum as for the collapsed sum of permutations to a sum over indices. Hence, starting with the last factor, for $i_1\neq \cdots \neq i_{m-1}$ the sum equals,
\begin{align}
 \sum_{i_1\neq ... \neq i_m}  \hat{x}_{i_1}^{n_1}...\hat{x}_{i_m}^{n_m}=&\\
 =& \sum_{i_1\neq ... \neq i_{m-1}} \hat{x}_{i_1}^{n_1}...\hat{x}_{i_{m-1}}^{n_{m-1}}\sum_{i_m \neq i_1,...,i_{m-1}} \hat{x}_{i_m}^{n_m} \nonumber\\
 =&\sum_{i_1\neq ... \neq i_{m-1}} \hat{x}_{i_1}^{n_1}...\hat{x}_{i_m}^{n_{m_1}} \Big [ n\langle \hat{x}^{n_m}\rangle - \hat{x}_{i_1}^{n_m} - ... - \hat{x}_{i_{m-1}}^{n_m}    \Big ] \nonumber
\end{align}
Proceeding in this manner gives us the inductive formula for the $x$ terms,
\begin{equation}
    X_{(n_1,..,n_m)}^{m,k}  =   \ n \langle  \hat{x}^{n_m} \rangle X_{(n_1,...,n_{m-1})}^{m-1,k-n_m}  - \sum_j X_{(n_1 + n_m\delta_{ij},...,n_{m-1}+n_m\delta_{(m-1)j})}^{m-1,k} 
    \label{inducx2}
\end{equation}
The adjusted multinomial coefficient represents the number of ways to place $k$ distinct object into $m$ distinct groups, of respective size $n_1,\cdots,n_m$, with the caveat that groups of the same size are \textit{indistinguishable}. Thus, partition the exponents $(n_1,\dots,n_m)$ into subgroups, $g_1,\dots,g_r$, via the equivalence relation $n_i\equiv n_j \iff n_i=n_j$, and let $d_i=|g_i|$. The over count that arises, permuting each group, is equal to $d_1!\cdots d_r!$, so that 
\begin{equation}
    \binom{k}{n_1,\dots,n_m}^* = \ \frac{1}{d_1!}\cdots\frac{1}{d_r!}\binom{k}{n_1,\dots,n_m}
\end{equation}

To understand why this adjustment is needed, we proceed with an example. Let $k=4$, $m=2$, and $n_1=n_2=2$ and look at the ways the terms in 
\begin{equation}
    \sum_{i \neq j} \hat{x}_{i}^2\hat{x}_{j}^2
    \label{ex22}
\end{equation}
arise from the original product of sums on the left side of \eqref{xterms}. For fixed $i\neq j$, We have to pull two $\hat{x}_{i}$ and two $\hat{x}_{j}$  from the product of the four distinct sums
\begin{equation}
   \Big [ \sum_l \hat{x}_{l} \Big ] \Big [ \sum_r \hat{x}_{r} \Big ] \Big [ \sum_s \hat{x}_{s} \Big ] \Big [ \sum_t \hat{x}_{t} \Big ]
\end{equation}
for which there are $6 = \binom{4}{2,2}$ ways to do this
\begin{eqnarray}
    \hat{x}_{i}  \hat{x}_{i} \hat{x}_{j}  \hat{x}_{j} \ \ \ \hat{x}_{j}  \hat{x}_{j} \hat{x}_{i}  \hat{x}_{i} \nonumber\\
        \hat{x}_{i}  \hat{x}_{j} \hat{x}_{i}  \hat{x}_{j} \ \ \ \hat{x}_{j}  \hat{x}_{i} \hat{x}_{j}  \hat{x}_{i} 
          \label{ex22p}\\
          \hat{x}_{i}  \hat{x}_{j} \hat{x}_{j}  \hat{x}_{i} \ \ \ \hat{x}_{j}  \hat{x}_{i} \hat{x}_{i}  \hat{x}_{j} \nonumber
\end{eqnarray}
However, in \eqref{ex22p} to form a complete sum each row collapses into a single complete sum, for a total of $3= \frac{1}{2!}\binom{4}{2,2}$ distinct sums.

Putting everything together gives the desired induction formula for the moments of Pearson's distribution over $\Pi(D_n)$.
\begin{center}
 $\qed$   
\end{center}
\noindent
For $k=1$, we have 
\begin{equation}
\langle r_{\Pi} \rangle = \frac{n^2 (n-1)!\langle \hat{x} \rangle \langle \hat{y} \rangle } {n\cdot n!\hat{\sigma}_x \hat{\sigma}_y} ,
 \end{equation}
where by definition $\langle \hat{x}\rangle =\langle \hat{y} \rangle = 0$ , so that $\langle r_{\Pi} \rangle=0$ and consequently $\langle r_{\Pi}^k \rangle$ are also central moments. 
\\ 
\\
See Appendix \ref{app:form} for the analytic expansions for the first five moments of $r_{\Pi}$ and numerical validation. Additionally, see Appendix \ref{app:coef} for a method organizing the analytic expression for each moment.
 
%
% % % % % % % % % % % % % % % % % % % % % % % % % 

% % % % % % % % % % % % % % % % % % % % % % % % %
% Validation

% % % % % % % % % % % % % % % % % % % % % % % % % 
%End Notes

% % % % % % % % % % % % % % % % % % % % % % % % % 
% Conclusion
\section{Conclusion} 
Even without placing these formulas into a statistical framework, it is clear that a permutation test can not be distribution-free, as some methods for computing p-value assume. The permutation test depends on the central moments of the data, as previously investigated in \cite{hayes1996permutation}. The exact dependency of the permutation distribution on data and implications for asymptotic approximations as sample size increases, can potentially now be investigated \cite{aheizer1965classical,john2007techniques,mnatsakanov2008hausdorff}. The problem of approximating a distribution function given a fixed set of moments has been explored fairly thoroughly\cite{john2007techniques, mnatsakanov2008hausdorff}. Given a sufficient number of sample moments for the sample correlation over permutations of the data, one can estimate the distribution of the sample correlation over the set of permutations. Such an analysis could potentially be more efficient with the right implementation, certainly more efficient than computing permutations of data and their correlation coefficient. We posit that such a method exists for determining a direct p-value estimate.

%%%%%%%%%%%%%%%%%%%%%%%%%%%%%%%%%%%%%%%%%%%%%%
% Acknowledgements
%%%%%%%%%%%%%%%%%%%%%%%%%%%%%%%%%%%%%%%%%%%%%%

%%%%%%%%%%%%%%%%%%%%%%%%%%%%%%%%%%%%%%%%%%%%%%
% Declaration of interest
%%%%%%%%%%%%%%%%%%%%%%%%%%%%%%%%%%%%%%%%%%%%%%
%\section{Declaration of Interest}
%No potential conflict of interest is reported by the authors.

%%%%%%%%%%%%%%%%%%%%%%%%%%%%%%%%%%%%%%%%%%%%%%
% Bibliography
%%%%%%%%%%%%%%%%%%%%%%%%%%%%%%%%%%%%%%%%%%%%%%
\bibliographystyle{ieeetr}
\bibliography{main.bib}

%%%%%%%%%%%%%%%%%%%%%%%%%%%%%%%%%%%%%%%%%%%%%%
%% Appendix
%%%%%%%%%%%%%%%%%%%%%%%%%%%%%%%%%%%%%%%%%%%%%%
\begin{appendix}
\newpage
\section{Exact Formulas for $k=1,..,5$}
\label{app:form}
\noindent
From the induction formula every sample moment can be analytically determined. Below are the exact formulas for the first five sample moments in terms of the easily computed moments of the data. As in the main body of the paper, $\hat{\sigma}_x$ and $\hat{\sigma}_y$ are defined by \eqref{sigma}. Additionally, we employ the following notation for simplicity sake with $\chi_k=\langle x^k\rangle$, $\nu_k=\langle y^k \rangle$, and $\mu_{k,j}=\langle x^k \rangle \langle y^j \rangle$, all central moments. The first five sample moments of the sample correlation are as follows:
\begin{eqnarray*}
    \langle r_{\Pi}^1 \rangle  = &\  0
    \\
    \langle r_{\Pi}^2 \rangle = &\ \frac{1}{(n-1)}\\
    \\
     \langle r_{\Pi}^3 \rangle  = &\  \frac{\mu_{3,3}}{\hat{\sigma}_x^3\hat{\sigma}_y^3}\bigg[\frac{1}{n^2}h_{n,1} +\frac{3}{n^2(n-1)}h_{n,2} +\frac{4}{n^2(n-1)(n-2)}h_{n,3}\bigg]
\\
\\
      \langle r_{\Pi}^4 \rangle  = & \ \frac{1}{\hat{\sigma}_x^4 \hat{\sigma}_y^4} \bigg[\frac{\mu_{4,4}}{n^3}h_{n,1} + 
        \frac{4\chi_4 \nu_4}{n^3(n-1)}h_{n,2} +
          \frac{3[n^2\hat{\sigma}_x^4-n\chi_4][n^2\hat{\sigma}_y^4 -n\nu_4]}{n^5(n-1)}h_{n,2}
        \\ \\
        & +\frac{6[2n\chi_4 -n^2\hat{\sigma}_x^4][2n\nu_4 -n^2\hat{\sigma}_y^4]}{n^5(n-1)(n-2)}h_{n,3}
        + \frac{9[2n\chi_4 - n^2 \hat{\sigma}_x^4][2n\nu_4 - n^2 \hat{\sigma}_y^4]}{n^5(n-1)(n-2)(n-3)}h_{n,4}\bigg] 
\end{eqnarray*}
\begin{eqnarray*} 
       \langle r_{\Pi}^5 \rangle  = & \ \frac{1}{\hat{\sigma}_x^5 \hat{\sigma}_y^5} \bigg[ \frac{\mu_{5,5}}{n^4}h_{n,1} 
         + 5\frac{\mu_{5,5}}{n^4(n-1)}h_{n,2}
         +
         10\frac{[n^2\chi_3 \chi_2 -n\chi_5][n^2\nu_3 \nu_2 -n\nu_5]}
         {n^6(n-1)}h_{n,2}
         \\ 
         \\ 
         &
         +10\frac{[2n\chi_5 -n^2\chi_3\chi_2][2n\nu_5 -n^2\nu_3 \nu_2]}{n^6(n-1)(n-2)}h_{n,3}
         \\
         \\
         & + 60\frac{[n\chi_5-n^2\chi_3\chi_2][n\nu_5-n^2 \nu_3 \nu_2]}{n^6(n-1)(n-2)} h_{n,3}\\
         \\
         & +10\frac{[6n\chi_5 - 5n^2\chi_3\chi_2][6n\nu_5 - 5n^2\nu_3\nu_2]}{n^6(n-1)(n-2)(n-3)}h_{n,4
         }\\
         \\& +
         \frac{16[6n\chi_5 - 5n^2\chi_3\chi_2][6n\nu_5 - 5n^2\nu_3\nu_2]}{n^6(n-1)(n-2)(n-3)(n-4)}h_{n,5}\bigg]
\end{eqnarray*}

%\subsection{Validation}
%\newpage
\subsection{Validation}
\label{app:valid}
In order to demonstrate the validity of these formulas, we randomly generated datasets of sizes $n=\{3,\dots,8\}$ over $100$ trials at each $n$. We compared the derived moments to the empirical moments computed directly from $\Pi(D_n)$, which can be done in the case of small datasets, however for larger datasets ($n>8$) this becomes computationally impractical.
\\
\\
The mean squared error was computed as:
\begin{equation}
  MSE = \frac{1}{N_{trials}}\sum \left(\langle r^k\rangle_{\Pi(D_n)} - \langle r^k\rangle_{exact}\right)^{2}
\end{equation}
\\
This validation procedure was performed using MATLAB with double-precision floating point computations. The errors tabulated in Table \ref{tab:validation} are within machine epsilon error indicating that the formulas are indeed exact.

\begin{table}[h]
\centering
\begin{tabular}{| c || c | c | c | c | }

\hline
 &\multicolumn{4}{|c|}{MSE of $k^{th}$ Moment } \\ [0.5ex] 
 
  Sample Size & 2&3&4&5 \\ [0.5ex] 
 
 \hline
 $n=3$ & $1.20\mathrm{e}{-32}$ & $2.57\mathrm{e}{-32}$ & $3.99\mathrm{e}{-32}$ & $4.82\mathrm{e}{-32}$\\ 
 \hline
 $n=4$ & $6.81\mathrm{e}{-33}$ & $2.46\mathrm{e}{-33}$ & $1.11\mathrm{e}{-32}$ & $3.18\mathrm{e}{-33}$\\ 
 \hline
 $n=5$ & $6.06\mathrm{e}{-33}$ & $1.23\mathrm{e}{-33}$ & $5.75\mathrm{e}{-32}$ & $1.55\mathrm{e}{-33}$\\ 
 \hline
 $n=6$ & $2.42\mathrm{e}{-32}$ & $8.99\mathrm{e}{-34}$ & $6.00\mathrm{e}{-33}$ & $4.37\mathrm{e}{-34}$\\
 \hline
 $n=7$ & $1.31\mathrm{e}{-31}$ & $3.50\mathrm{e}{-33}$ & $1.78\mathrm{e}{-32}$ & $1.29\mathrm{e}{-33}$\\
 \hline
 $n=8$ & $7.38\mathrm{e}{-31}$ & $1.05\mathrm{e}{-32}$ & $5.37\mathrm{e}{-32}$ & $3.17\mathrm{e}{-33}$\\
 \hline
 \end{tabular}
\caption{Validation Error}
\label{tab:validation}
\end{table}

%\newpage
\section{Organizing the computation of $\langle r_{\Pi}^5\rangle$}
\label{app:coef}
\noindent
Assuming the number of data points $n>5$ and so ignoring the $h_{n,m}$ factor, from \eqref{pcpm} and for $k=5$ we have:
\begin{align}
    \langle r_{\Pi}^5 \rangle  & = &  \nonumber \\ 
    & = & \frac{1}{n^k n!\hat{\sigma}_x^k \hat{\sigma}_y^k} \sum_{m=1}^{5}  \ \ \sum_{n_1+...+n_m=5} \binom{5}{n_1,..,n_m}^* X_{(n_1,..,n_m)}^{m,5} (n-m)!Y_{(n_1,..,n_m)}^{m,5}
    \label{apxb}
\end{align}
The sum is organized as follows:
\begin{align*}
\begin{tikzpicture}[scale=.7]
    % Center nodes
    \node (lc) at (0,9.5) {$Z_\mathbf{(n_1,\dots,n_m)}^{m,k}$};
    \node (11111) at (0,8) {$Z_{(1,1,1,1,1)}^{5,5}$};
    \node (2111) at (0,6)  {$Z_{(2,1,1,1)}^{4,5}$};
    \node (221) at (-2,4)  {$Z_{(2,2,1)}^{3,5}$};
    \node (311) at (2,4) {$Z_{(3,1,1)}^{3,5}$};
    \node (32) at (-2,2) {$Z_{(3,2)}^{2,5}$};
    \node (41) at (2,2)  {$Z_{(4,1)}^{2,5}$};
    \node (5) at (0,0)  {$Z_{(5)}^{1,5}$};
    \draw [->] (11111) -- (2111);
    \draw [->] (2111) -- (221);
    \draw [->] (2111) -- (311);
    \draw [->] (221) -- (32);
    \draw [->] (221) -- (41);
    \draw [->] (311) -- (32);
    \draw [->] (311) -- (41);
    \draw [->] (41) -- (5);
    \draw [->] (32) -- (5);
    % Right labels
    \node (r1) at (6,9.5) {$\mathbf{\Pi}$-multiplicity};
    \node (r1) at (6,8) {$(n-5)!$};
    \node (r1) at (6,6) {$(n-4)!$};
    \node (r1) at (6,4) {$(n-3)!$};
    \node (r1) at (6,2) {$(n-2)!$};
    \node (r1) at (6,0) {$(n-1)!$};
    % Left labels
   \node (l1) at (-6,9.5) {$\mathbf{\binom{k}{n_1,..,n_m}^*}$};
    \node (l1) at (-6,8) {$1$};
    \node (l1) at (-6,6) {$10$};
    \node (l1) at (-6,4) {$15/10$};
    \node (l1) at (-6,2) {$10/5$};
    \node (l1) at (-6,0) {$1$};
\end{tikzpicture}
\end{align*}
\\
For $Z=X$ or $Y$ we have,
\begin{eqnarray}
Z_{(1,1,1,1,1)}^{5,5} & = & \sum \hat{z}_{i_1} \hat{z}_{i_2} \hat{z}_{i_3} \hat{z}_{i_4} \hat{z}_{i_5} \nonumber \\
  & = & \sum \hat{z}_{i_1} \hat{z}_{i_2} \hat{z}_{i_3} \hat{z}_{i_4}(n\langle z \rangle - \hat{z}_{i_1} - \hat{z}_{i_2}- \hat{z}_{i_3} -\hat{z}_{i_4}) \nonumber \\
  & = & -4Z_{(2,1,1,1)}^{4,5} \nonumber \\
  %%%%%%%%%%%%%%%%%%%%%%%%%%%%%%%%%%%%%%%%%%%%%%
Z_{(2,1,1,1)}^{4,1} & = & \sum \hat{z}_{i_1}^2 \hat{z}_{i_2} \hat{z}_{i_3} \hat{z}_{i_4}  \nonumber \\
  & = & \sum \hat{z}_{i_1}^2 \hat{z}_{i_2} \hat{z}_{i_3} (n\langle z \rangle - \hat{z}_{i_1} - \hat{z}_{i_2}- \hat{z}_{i_3}) \nonumber \\
  & = & -Z_{(3,1,1)}^{3,5}  - 2Z_{(2,2,1)}^{3,5}\nonumber \\
  %%%%%%%%%%%%%%%%%%%%%%%%%%%%%%%%%%%%%%%%%%%%%%%%%%%%%%%
Z_{(2,2,1)}^{3,5} & = & \sum \hat{z}_{i_1}^2 \hat{z}_{i_2}^2 \hat{z}_{i_3}  \nonumber \\
  & = & \sum \hat{z}_{i_1}^2 \hat{z}_{i_2}^2 \hat{z}_{i_3} (n\langle z \rangle - \hat{z}_{i_1} - \hat{z}_{i_2}) \nonumber \\
  & = & -2Z_{(3,2)}^{2,5}\nonumber \\
%%%%%%%%%%%%%%%%%%%%%%%%%%%%%%%%%%%%%%%%%%%%%%%%%%%%%%%
Z_{(3,1,1)}^{3,5} & = & \sum \hat{z}_{i_1}^3 \hat{z}_{i_2} \hat{z}_{i_3}  \nonumber \\
  & = & \sum \hat{z}_{i_1}^3 \hat{z}_{i_2} \hat{z}_{i_3} (n\langle z \rangle - \hat{z}_{i_1} - \hat{z}_{i_2}) \nonumber \\
  & = & -Z_{(4,1)}^{2,5}  - Z_{(3,2)}^{2,5}\nonumber \\
%%%%%%%%%%%%%%%%%%%%%%%%%%%%%%%%%%%%%%%%%%%%%%%%%%%%%%%
Z_{(3,2)}^{2,5} & = & \sum \hat{z}_{i_1}^3 \hat{z}_{i_2}^2  \nonumber \\
  & = & \sum \hat{z}_{i_1}^3  (n\langle z^2 \rangle - \hat{z}_{i_1}^2 ) \nonumber \\
  & = & n^2\langle z^2 \rangle Z_{(3)}^{1,3}   - Z_{(5)}^{1,5}\nonumber \\
  %%%%%%%%%%%%%%%%%%%%%%%%%%%%%%%%%%%%%%%%%%%%%%%%%%%%%%%
Z_{(4,1)}^{2,5} & = & \sum \hat{z}_{i_1}^4 \hat{z}_{i_2}  \nonumber \\
  & = & \sum \hat{z}_{i_1}^4  (n\langle z \rangle - \hat{z}_{i_1} ) \nonumber \\
  & = &  -Z_{(5)}^{1,5}\nonumber \\
  %%%%%%%%%%%%%%%%%%%%%%%%%%%%%%%%%%%%%%%%%%%%%%%%%%%%%%%
  Z_{(5)}^{1,5} & = & \sum \hat{z}_{i_1}^5 \nonumber \\ 
    & = & n\langle z^5\rangle \nonumber \text{, where in general } Z_{(m)}^{1,m} = n\langle z^m\rangle
\end{eqnarray}
Combining the terms we have,
\begin{eqnarray}
\langle r_{\Pi}^5 \rangle %& = &\nonumber \\
& = & \frac{1}{n^k n! \hat{\sigma}_x^k \hat{\sigma}_y^k}\Big[ (n-1)!X_{(5)}^{1,5} Y_{(5)}^{1,5} \nonumber \\ & \ & + 5(n-2)!X_{(4,1)}^{2,5}Y_{(4,1)}^{2,5} + 10(n-2)!X_{(3,2)}^{2,5}Y_{(3,2)}^{2,5}  \nonumber \\
& \ &  + 10(n-3)!X_{(3,1,1)}^{3,5}Y_{(3,1,1)}^{3,5} + 15(n-3)!X_{(2,2,1)}^{3,5}Y_{(2,2,1)}^{3,5}\nonumber \\
& \ & + 10(n-4)!X_{(2,1,1,1)}^{4,5}Y_{(2,1,1,1)}^{4,5} + (n-5)!X_{(1,1,1,1,1)}^{5,5}Y_{(1,1,1,1,1)}^{5,5}\Big ] \nonumber
\end{eqnarray}
\noindent
From which it follows, as seen in Appendix \ref{app:form}, that:
\begin{eqnarray*} 
       \langle r_{\Pi}^5 \rangle  = & \ \frac{1}{\hat{\sigma}_x^5 \hat{\sigma}_y^5} \bigg[ \frac{\mu_{5,5}}{n^4}
         + 5\frac{\mu_{5,5}}{n^4(n-1)}
         +
         10\frac{[n^2\chi_3 \chi_2 -n\chi_5][n^2\nu_3 \nu_2 -n\nu_5]}
         {n^6(n-1)}
         \\ 
         \\ 
         &
         +10\frac{[2n\chi_5 -n^2\chi_3\chi_2][2n\nu_5 -n^2\nu_3 \nu_2]}{n^6(n-1)(n-2)}
         + 60\frac{[n\chi_5-n^2\chi_3\chi_2][n\nu_5-n^2 \nu_3 \nu_2]}{n^6(n-1)(n-2)} \\
         \\
         & +10\frac{[6n\chi_5 - 5n^2\chi_3\chi_2][6n\nu_5 - 5n^2\nu_3\nu_2]}{n^6(n-1)(n-2)(n-3)}\\
         \\& +
         \frac{16[6n\chi_5 - 5n^2\chi_3\chi_2][6n\nu_5 - 5n^2\nu_3\nu_2]}{n^6(n-1)(n-2)(n-3)(n-4)}\bigg]
\end{eqnarray*}
%%\section{\texorpdfstring{$(1,1,1,1,1)$}{1,1,1,1,1}}

\end{appendix}
% https://statistics.stanford.edu/sites/g/files/sbiybj6031/f/2015-15.pdf

\end{document}